# Computation and applications of limits of certain non-stationary Markov chains

Samuel Awoniyi
Department of Industrial and Manufacturing Engineering
FAMU-FSU College of Engineering
Tallahassee, Florida
ORCID# 0000 0001 7102 6257

# Abstract

This article presents an algorithm for computing limits of a class of non-stationary Markov chains. Interest in this class was motivated by a need to compute certain healthcare cycles. A mathematical validation of the algorithm is given. Application instances are described that include predicting healthcare cycles in a nutritional healthcare center, and predicting cycles of a general system maintenance. To be able to apply the algorithm, one only has to be able to put treatment goals into a {0,1} matrix that we refer to as tg-matrix. An objective of this article is to hopefully foster investigation and applications of various practical non-stationary Markov chains.

# 1. Introduction

In its general form, Markov chain is a mathematical model suited to describing and predicting aspects of various real-world non-deterministic processes, such as preventative healthcare cycles and machine/system maintenance cycles [1,4,6,8]. However, current literature on Markov chain is somewhat skewed to and concentrated on stationary Markov chains, that is, Markov chains under the stationarity assumption. The stationarity assumption states that initial Markov chain data (that is, sojourn time at each state and one-step transition probabilities at each state) do not change over the time period of interest.

Perhaps one explanation for that literature bias is that computations for stationary Markov chains are relatively straightforward when compared to computations for non-stationary Markov chains. Furthermore, there are widely disparate ways in which non-stationary Markov chains may be described [3,5,7,10].

In real-world Markov chain modeling wherein underlying non-deterministic processes can learn, adjust or heal over time, as in biological and healthcare-related processes, one should not lightly make the stationarity assumption, even though the stationarity assumption could simplify requisite analysis or computations. Instead, one should consider using a simulation method to approximate and analyze the nature of non-stationarity involved, and thereby obtain a possibly useful limit of underlying non-stationary Markov chain. This article is a contribution made in that vein.

This article describes an algorithm for computing limits of non-stationary continuous-time Markov chains suited to modeling preventative healthcare cycles and similar general system maintenance cycles. This class of non-stationary Markov chains is one wherein each state sojourn time does not change during the time period of interest, but the associated one-step transition probabilities change uniformly during the time period of interest. This article also describes briefly two broad application instances, namely, a nutritional healthcare center (a prototype product-tracking example), and a general system maintenance instance, a prototype system-tracking example.

The remainder of this article is organized as follows. Section 2 introduces a new class of stationary Markov chains (MC's) and its limit characteristics. Section 3 describes a class of non-stationary MC's whose limit is the main subject matter of this article. Section 4 describes an algorithm for computing limits of that class of non-stationary MC's. Section 5 describes two broad application instances. Section 6 provides a rigorous validation of the algorithm stated in Section 4.

# 2. A class of stationary MCs

In this Section we introduce a class of stationary continuous-time MC's which we find to be the limit of the class of non-stationary MC's that is the main subject of this article. We will refer to this class of stationary continuous-time MC's as "treatment goal continuous-time MC's" (to be abbreviated simply as "tg MC"). We will also briefly discuss the nature of associated sojourn-time cycles.

In this article, any n-state non-stationary continuous-time Markov chain will be represented as an infinite sequence, $\{[T_{(t)}, P_{(t)}],\ t = 1, 2, \ldots\}$, of stationary continuous-time Markov chains, wherein the ordered pair $[T_{(t)}, P_{(t)}]$ is an n-state stationary continuous-time Markov chain having n-vector $T_{(t)}$ of (state-by-state) mean sojourn-times, and having n-by-n matrix, $P_{(t)}$, of embedded MC's one-step transition probabilities. In the remainder of this article, $T_{(t)}$ is a constant vector, and we will accordingly write it as $T$; we will sometimess write $P_{(t)}$ simply as $P$, in the interest of notation tidiness.

**Definition** *A stationary continuous-time MC in sojourn-time form, say [T,P], is called a "treatment goal MC" if each row of the stochastic matrix P is a unit n-vector (that is, it has 1 in exactly one component and 0 in every other component). The matrix P will sometimes be refered to as a tg-matrix.*

An example of such a tg-matrix is

$$P = \begin{pmatrix} 0 & 0 & 1 & 0 & 0 \\ 0 & 0 & 0 & 0 & 1 \\ 0 & 0 & 0 & 1 & 0 \\ 1 & 0 & 0 & 0 & 0 \\ 0 & 1 & 0 & 0 & 0 \end{pmatrix}$$

Every element in the diagonal of $P$ is 0, because $[T, P]$ is a continuous-time MC.

In this article, we will let emb($[T, P]$) denote the embedded MC of $[T, P]$; note that emb($[T, P]$) is a MC that is completely described by $P$. For the particular instance of $P$ stated above, the (system of) balance equations of emb($[T, P]$) is, with *tr* meaning "matrix transpose",

$$\left\{ \begin{array}{l} P^{tr}x = x \\ x_1 + x_2 + x_3 + x_4 + x_5 = 1 \\ x_i \geq 0, i = 1, .., 5 \end{array} \right\} \cdots\cdots(Eq)$$

That particular 5-state instance system of balance equations happens to have two solutions, $x = (1/3,\ 0,\ 1/3,\ 1/3,\ 0)^{tr}$ and $x = (0,\ 1/2,\ 0,\ 0,\ 1/2)^{tr}$.

As illustrated in the network diagram below, the solution $x = (1/3,\ 0,\ 1/3,\ 1/3,\ 0)^{tr}$ corresponds to the cycle $1 \to 3 \to 4 \to 1$ in the network of emb($[T, P]$), and the solution $x = (0,\ 1/2,\ 0,\ 0,\ 1/2)^{tr}$ corresponds to the cycle $2 \to 5 \to 2$.

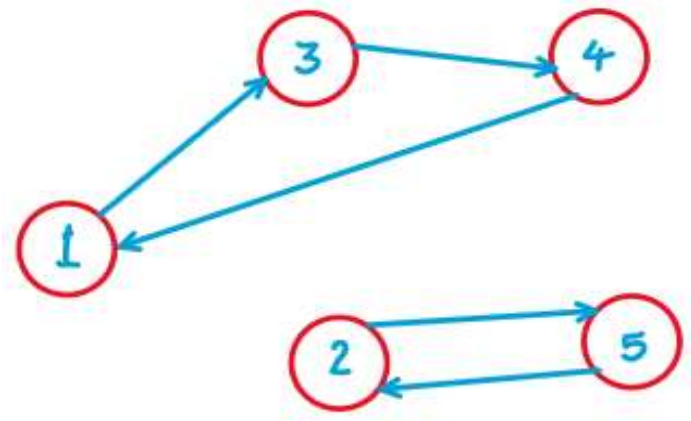

If it starts at any state in $\{1,3,4\}$, then the transitions of that $[T,P]$ will always follow the cycle $1 \rightarrow 3 \rightarrow 4 \rightarrow 1$. With T=$(t_1,t_2,t_3,t_4,t_5)^{tr}$, one could trace that network walk of $[T,P]$ to compute $[T,P]$'s mean sojourn times in any one of the subsets $\{1\}$, $\{3\}$, $\{4\}$, $\{1,3\}$, $\{1,4\}$, $\{3,4\}$. One would then get the following table, wherein s(G) denotes sojourn time inside G.

| G | $\{1\}$ | $\{3\}$ | $\{4\}$ | $\{1,3\}$ | $\{1,4\}$ | $\{3,4\}$ |
|---|---|---|---|---|---|---|
| s(G) | $t_1$ | $t_3$ | $t_4$ | $t_1+t_3$ | $t_4+t_1$ | $t_3+t_4$ |

Note that, for that table, the sojourn-time cycle for any subset G, s(G)+s($G^c$), is $t_1 + t_3 + t_4$. Similarly, if it starts at any state in $\{2,5\}$, the transitions of $[T,P]$ will always follow the cycle $2 \rightarrow 5 \rightarrow 2$, and one can follow that $[T,P]$ walk to compute its mean sojourn times in any one of the subsets $\{2\}$, $\{5\}$.

In general, the balance equations of a tg MC can have multiple solutions, each one corresponding to a cycle in a network representation of the tg MC. The number of solutions of such a system of balance equations is the number of cycles contained in the network representation of the tg MC, and all the solutions of such a system of balance equations may be efficiently computed with a simple generalization of the GTH algorithm [2]. A MATLAB implementing such a generalization of the GTH algorithm is included in this article as an Appendix.

# 3. A class of non-stationary MCs

We describe in this Section a class of non-stationary continuous-time MC's whose limits turn out to be tg MC's. This class of non-stationary MC's is motivated by a need to compute certain healthcare cycles. First, we give an informal description. Thereafter, we will utilize a discrete-event simulation to give a formal description.

## 3.1 Informal description

We describe here a class of non-stationary MCs, say $\{[T,P_{(t)}],\ t = 1,2\ldots\}$, with the following three defining properties:

(i) $T$, a sojourn time vector, does not change during the time period of interest, that is, for all values of $t$;

(ii) $P_{(t)}$, a one-step transition matrix of emb($[T,P_{(t)}]$) "*changes uniformly in consonance with the starting value of P, say* $P_{(1)}$", beginning at the starting row, row $s_{(1)}$, of $P_{(1)}$;

(iii) each row of $P_{(1)}$ has exactly one component that contains the maximum probability for that row, and, or all large values of $t$, the component containing the largest probability in each row of $P_{(t)}$ remains the same as it is in $P_{(1)}$,

Regarding property (ii), we explain what is meant by "changes uniformly in consonance with an initial value, say $P_{(1)}$" as follows. A discrete-event simulation run, utilizing a uniform distribution random variate, generates a sequence $\{[s_{(t)},\ P_{(t)}],\ t = 1,2\ldots\}$, with each $P_{(t)}$ a stochastic matrix, and $s_{(t)}$ an index of a row of $P_{(t)}$. The matrix $P_{(t)}$ gets transformed into matrix $P_{(t+1)}$ through a change in row $s_{(t)}$ only. The following table describes an iteration of that discrete-event simulation.

| *step 0* | *step 1* | *step 2* | *step 3* |
|---|---|---|---|
| start with $[s_{(t)}, P_{(t)}]$ as inputs | pick random number $r \in [0,1]$ | inside row $s_{(t)}$ of $P_{(t)}$, $r$ points to row index $s_{(t+1)}$ | adjust row $s_{(t)}$ of $P_{(t)}$, thereby obtaining $P_{(t+1)}$ |

It turns out that the sequence $\{P_{(t)},\ t = 1,2\ldots\}$ converges to a tg-matrix, and that largest probability (in each row) converges to 1.

In the healthcare cycle application instance, the maximal component of each row of $P_{(1)}$ indicates the *desirable outcome* or the *expected effect of healing and learning* over time, as a consequence of planned health-related actions.

## 3.2 Formal description

We present here a quantitative version of the qualitative description presented above, along with an instructive numerical illustration.

The non-stationary Markov chain whose limit is the subject of this article is formally denoted here by the sequence $\{[T, P_{(t)}],\ t = 1,2,\ldots,\}$, wherein $T$ is a constant n-vector of mean-sojourn-times at the states, and $P_{(1)},\ P_{(2)},\ldots$ is a sequence of one-step transition probability matrices (of embeded MC's), as specified in the following several paragraphs.

The initial matrix $P_{(1)}$ is assumed to be given, along with a starting row index, say $s_{(1)}$, and a probability simulation change parameter, say ε, is chosen such that $0 < \varepsilon < 1$.

For $t = 1,2,\ \ldots,$ the matrix $P_{(t)}$ changes into the matrix $P_{(t+1)}$ through a "randomly generated" change in one row only, the row with index $s_{(t)}$. That randomly generated change produces the next row index $s_{(t+1)}$ as well. The random change in row $s_{(t)}$ is effected through the following sub-procedure.

*Sub-procedure for effecting changes in row* $s_{(t)}$ : With $P_{(t)}$ and $s_{(t)}$ as input, *first* pick a pseudo-random number $r \in [0,1]$, and then use $r$ to "point-to" position $(s_{(t)},j)$ in the matrix $P_{(t)}$, as explained in "Definition of point-to" below; then set $s_{(t+1)} \leftarrow j$. *Secondly*, multiply the current probability number in position $(s_{(t)},s_{(t+1)})$ by $1 + \varepsilon$; thereafter, normalize the resultant row $s_{(t)}$ of $P_{(t)}$, by dividing (that) row $s_{(t)}$ by the sum of all its elements. *Finally*, set $P_{(t+1)} \leftarrow$"the changed matrix $P_{(t)}$". Thus, obtaining each term of the sequence $P_{(2)},\ P_{(3)},\ldots$ may be regarded as the result of an event of a discrete-event simulation run that utilizes a U(0,1) random variate for each row change.

*Definition of point-to:* A pseudo-random number value $r$ is said to "point-to" position $(s_{(t)},j)$ in row $s_{(t)}$ of matrix $P_{(t)}$ if $\sum_{k=1}^{k=j-1} P_{(t)}(s_{(t)},k) < r \leq \sum_{k=1}^{k=j} P_{(t)}(s_{(t)},k)$.

A less formal way of stating that definition of "point-to" is to say that elements of row $s_{(t)}$ of matrix $P_{(t)}$ are summed up, from left to right, until the sum is no smaller than the value $r$, with that happenning for the first time in column $j$ of $P_{(t)}$. Note that the value $r$ will never point to diagonal position $(s_{(t)},s_{(t)})$, since each diagonal element of $P_{(t)}$ is 0, as $P_{(t)}$ is the one-step transition matrix for the embedded MC of a continuous-time Markov chain.

*A summary of the formal description:*

(i) Thus, changing row $s_{(t)}$ of $P_{(t)}$ yields $P_{(t+1)}$ along with its own row index $s_{(t+1)}$ (with $s_{(t+1)} \neq s_{(t)}$); that is, $P_{(t+1)}$ differs from $P_{(t)}$ only in row $s_{(t)}$.

(ii) Since $P_{(1)},\ P_{(2)},\ldots$ and $s_{(1)},s_{(2)},\ldots$ are sequences having infinitely many terms, whereas $P_{(1)}$ has only a finite number of rows, at least one row index will be repeated infinitely many times in the sequence of row indices $s_{(1)},s_{(2)},\ldots.$ That is, the sequence $s_{(1)},s_{(2)},\ldots$ contains at least one cycle.

(iii) The sequence $s_{(1)},s_{(2)},\ldots$ induces a sequence of minimal cycles, each one corresponding to

a cycle in the network of $P_{(1)}$. As an illustration, the sub-sequence 1, 2, 4, 1, 3, 4, 1, 4, 2 induces the minimal cycles (1, 2, 4, 1), (1, 3, 4, 1), (4, 1, 3, 4), (4, 1, 4); note that (2, 4, 1, 3, 4, 1, 4, 2) is a non-minimal cycle because it contains another cycle (4, 1, 3, 4).

# 4. Resultant algorithm

We next state our algorithm for computing limits of the class of non-stationary MC's described above. It consists of computing a cycle in the network of $P_{(1)}$.

ALGORITHM: *Under the assumption that the maximum probability in each row of $P_{(1)}$ occurs in only one component of that row, the Steps of this algorithm are: Step 1 - start at the node (in the network of $P_{(1)}$) corresponding to state $s_{(1)}$; and Step 2 - trace the path guided by the maximum probability coming out of each node, until a node is reached that is already included on the path. The existence of such a cycle is guaranteed by the fact that the network of $P_{(1)}$ is loopless (since the matrix $P_{(1)}$ has only zeros on its diagonal).*

As an illustration of this algorithm, suppose $s_{(1)} = 2$, and suppose $P_{(1)}$ is the following 6-by-6 stochastic matrix, with the maximum probabilty emphasized in each row.

| . | 1 | 2 | 3 | 4 | 5 | 6 |
|---|---|---|---|---|---|---|
| 1 | 0 | 1/2 | 0 | 0 | 1/3 | 1/6 |
| 2 | 0 | 0 | 1/2 | 0 | 1/4 | 1/4 |
| 3 | 0 | 0 | 0 | 1/2 | 1/8 | 1/4 |
| 4 | 0 | 0 | 0 | 0 | 1/9 | 8/9 |
| 5 | 1/2 | 1/4 | 1/4 | 0 | 0 | 0 |
| 6 | 1/4 | 1/2 | 1/8 | 1/8 | 0 | 0 |

Starting at $s_{(1)} = 2$, the algorithm produces the cycle 2 → 3 → 4 → 6 → 2. The network for that result is shown in Figure 1 below, with the cycle 2 → 3 → 4 → 6 → 2 displayed in broken arrows.

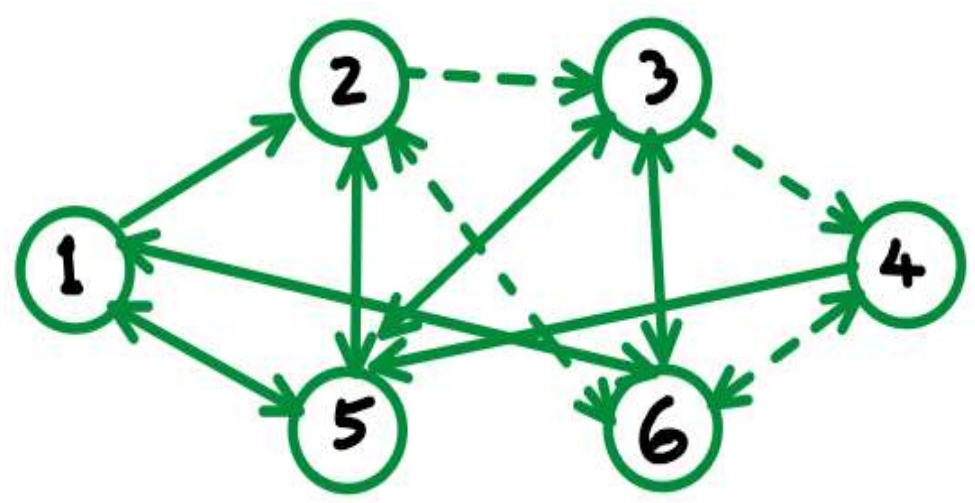

Figure 1

# 5. Two application instances

As an illustration of how to apply our algorithm described above, we will briefly describe in this section two broad application instances, namely, a nutritional healthcare center instance and a general system maintenance instance. In the nutritional healthcare center instance, the object being tracked by the MC is the center's product (the member's health outcome). In the general system maintenance instance, the object being tracked by the MC is the status of the production system, not the system's product.

## 5.1 Application 1 - A nutritional healthcare center instance

This is a hypothetical instance of a preventative nutritional healthcare center whose business essentially consists of specifying prescriptions (for food, medication, physical exercise schedules, associated home environment, etc) for center members/customers, with the goal of enabling good health for each customer as an individual.

### How the center operates

Persons become members of this nutritional healthcare center by registering and, thereafter, by completing every six months a set of relevant laboratory assessments (imaging, blood analysis, etc). Each time, the assessments are summarized into an initial data/matrix of a non-stationary Markov process for the purpose of possibly predicting aspects of the customer's health and prescribing suitable responses (food, medication, exercise, etc) as well.

During each six-month interval, center members perform minor assessments at home and take pertinent maintenance actions as outlined on the healthcare center's website. The efficacy of such maintenance actions is what translates into required Markov assumption, and the resulting learning and healing, on the part of each center member, is what translates into a Markov chain non-stationarity.

### An example of the questions addressed for each customer

For a particular customer, let us suppose that there are five possible "good states" G1, G2, G3, G4, G5, along with two possible "not-so-good states" B1, B2. Suppose that, for this customer, a center healthcare team has used a set of assessments to determine that the Markov chain sojourn times, denoted by a vector T, and the initial one-step transition probabilities (that is, $P_{(1)}$ as in Section 3), are as indicated in the following table.

.

| T | | . | G1 | G2 | G3 | G4 | G5 | B1 | B2 |
|---|---|---|---|---|---|---|---|---|---|
| 10 | ↔ | G1 | 0 | 1/2 | 1/6 | 0 | 0 | 1/6 | 1/6 |
| 12 | ↔ | G2 | 0 | 0 | 1/2 | 1/4 | 0 | 1/4 | 0 |
| 10 | ↔ | G3 | 0 | 0 | 0 | 3/7 | 2/7 | 1/7 | 1/7 |
| 14 | ↔ | G4 | 0 | 0 | 0 | 0 | 2/9 | 5/9 | 2/9 |
| 12 | ↔ | G5 | 0 | 0 | 0 | 0 | 0 | 7/9 | 2/9 |
| 1 | ↔ | B1 | 1/2 | 1/4 | 1/8 | 1/8 | 0 | 0 | 0 |
| 1 | ↔ | B2 | 1/8 | 1/2 | 1/8 | 1/8 | 1/8 | 0 | 0 |

.

One of the questions that may be addressed is: what is the *"the average length of time the customer may be expected to be in good states G1 or G2 or G3 or G4 or G5, compared to the average length of time the customer may be expected to be in not-so-good states B1or B2"* in the long run, given that the customer starts out from state G1?

### The algorithm's answer

It is clear, from the table displayed above, that the corresponding tg-matrix in this case, say $M$, is

.

$$M = \begin{pmatrix} 0 & 1 & 0 & 0 & 0 & 0 & 0 \\ 0 & 0 & 1 & 0 & 0 & 0 & 0 \\ 0 & 0 & 0 & 1 & 0 & 0 & 0 \\ 0 & 0 & 0 & 0 & 0 & 1 & 0 \\ 0 & 0 & 0 & 0 & 0 & 1 & 0 \\ 1 & 0 & 0 & 0 & 0 & 0 & 0 \\ 0 & 1 & 0 & 0 & 0 & 0 & 0 \end{pmatrix}$$

.

Applying our limit computation algorithm, we obtain the limit healthcare cycle $G1 \rightarrow G2 \rightarrow G3 \rightarrow G4 \rightarrow B1 \rightarrow G1$. We note how states $G5$ and $B2$ are not involved at the limit (that is, in the long run), a fact that is revealed only through our algorithm.

.

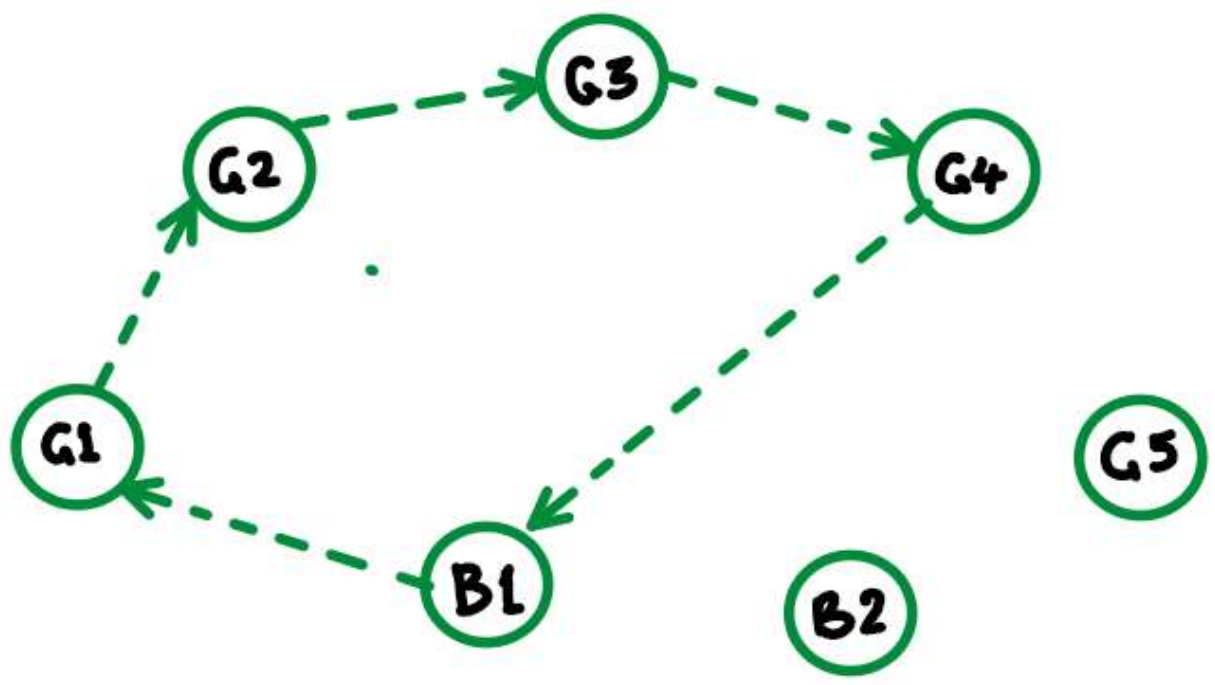

.

As notation towards answering the question posed above, let G ≡ {G1, G2, G3, G4, G5}, so that the complement $G^c$ ≡ {B1, B2}. Let s(G) denote the Markov chain's mean sojourn time inside G and let s($G^c$) denote the mean sojourn time inside $G^c$, as the Markov chain's transitions go back-and-forth between G and $G^c$ (see [1] to become familiar with this notation and some related computations).

Utilizing a computational method explained in [1], we then have s(G) = 46 and s($G^c$) = 1. Thus, the question posed above has been answered, along with the exact sequence of the Markov chain states involved..

**Remark** *For comparison purposes, suppose we made the assumption that the Markov chain data above remain valid for the customer for the foreseeable future, then we would utilize the procedure described in [1] to have s(G) = 24 and s($G^c$) = 1. The difference between s(G)=46 and s(G)=24 is a reflection of the long-term learning and healing on the part of the center customer.*

## 5.2 Application 2 - General system maintenance instance

We briefly describe here how to apply our algorithm to a general system maintenance instance. We consider a system that can be in one of $n$ states $S1, \ldots Sn$, with characteristics as follows:

– there is a mechanism for detecting/assessing system state change, whenever there is a state change;

– a state-dependent "treatment" is administered automatically as needed, for each state as soon as the state is detected;

– each state change is contingent upon the last-observed state only; this is to enable required Markov chain assumption;

– an n-by-n tg-matrix is available for the whole system.

In some general system application instances, on the basis of sojourn time information obtained for the non-stationary MC's limit, the system may be halted and adjusted in time before it enters into a very undesirable state.

# 6. A validation of the algorithm

We provide here a mathematical validation of the algorithm. This validation consists of explanations and summarizing lemmas that lead to the following concluding statement – if a row of $P_{(1)}$, say row $k$, is "visited" infinitely many times by the simulation run (described in Section 3, then the maximum probability in row $k$ converges to 1, thereby resulting in a row of a tg-matrix as defined in this article.

We will begin by using a small numerical example to illustrate what a row of $P_{(1)}$ becomes when its index is repeated infinitely many times in the sequence $s_{(1)}, s_{(2)}, \ldots$

## 6.1 On how rows of $P_{(i)}$ change

Here we utilize a general 4-by-4 instance of $P_{(1)}$ to analytically illustrate how rows of $P_{(t)}$'s are transformed when their indices are repeated infinitely many times in the sequence $s_{(1)}, s_{(2)}, \ldots$.

As definition, whenever there is a visit to a row, say row 1, in the sequence $s_{(1)}, s_{(2)}, \ldots$ (that is, whenever 1 occurs in the sequence $s_{(1)}, s_{(2)}, \ldots$), and a pseudo-random number thereafter points-to position $(1,j)$ of row 1, we shall then refer to position $(1,j)$ as *the pivot-position* for that row 1 visit.

We will use row 1 for this illustration. To that end, let us denote row 1 of $P_{(1)}$ by (0, $p_{1,2}$, $p_{1,3}$, $p_{1,4}$), and let us suppose that the first nine terms of $s_{(1)}, s_{(2)}, \ldots$ are 1, 2, 4, 1, 3, 4, 1, 4, 2. Accordingly, with $s_{(1)} = 1$, we have pivot positions (1,2), (1,3), (1,4), *in that order* (from the sub-sequence 1, 2, 4, 1, 3, 4, 1, 4, 2). Table 6.1 (below) displays what row 1 becomes immediately after each one of those first three visits to row 1.

Table 6.1

.

| | *what row 1 becomes immediately* |
|---|---|
| *after (1,2) as pivot-position* | $p_{1,1}$ is still 0<br>$p_{1,2}$ becomes $(1+\varepsilon)p_{1,2}/(1+\varepsilon p_{1,2})$<br>$p_{1,3}$ becomes $p_{1,3}/(1+\varepsilon p_{1,2})$<br>$p_{1,4}$ becomes $p_{1,4}/(1+\varepsilon p_{1,2})$ |
| *after (1,2), followed by (1,3) as pivot-positions* | $p_{1,1}$ is still 0<br>$p_{1,2}$ now becomes $(1+\varepsilon)p_{1,2}/(1+\varepsilon p_{1,2}+\varepsilon p_{1,3})$<br>$p_{1,3}$ now becomes $(1+\varepsilon)p_{1,3}/(1+\varepsilon p_{1,2}+\varepsilon p_{1,3})$<br>$p_{1,4}$ now becomes $p_{1,4}/(1+\varepsilon p_{1,2}+\varepsilon p_{1,3})$ |
| *after (1,2), followed by (1,3), followed by (1,4) as pivot-positions* | $p_{1,1}$ is still 0<br>$p_{1,2}$ becomes $p_{1,2}$ again<br>$p_{1,3}$ becomes $p_{1,3}$ again<br>$p_{1,4}$ becomes $p_{1,4}$ again |

Based on Section 3 description, the computation formula used to obtain Table 6.1 is:

.

"new value" = ["old value"]/[1+(ε. "current pivot-position value")],

.

except for current pivot-position itself, which becomes

.

[(1+ε). "current pivot-position value"]/[1+(ε. "current pivot-position value")].

.

For instance, $p_{1,2}$ in "*after (1,2), followed by (1,3) as pivot-positions"* is given by:

.

$$p_{1,2} \leftarrow [(1+\varepsilon)p_{1,2}/(1+\varepsilon p_{1,2})]/[1+\varepsilon p_{1,3}/(1+\varepsilon p_{1,2})]$$

.

which clearly reduces to

.

$$p_{1,2} \leftarrow (1+\varepsilon)p_{1,2}/(1+\varepsilon p_{1,2}+\varepsilon p_{1,3})$$

.

**Summary** *Note that the results in Table 6.1 will not change if we replace "after (1,2), followed by (1,3), as pivot-positions" with "after (1,3), followed by (1,2), as pivot-positions"; that is, when a row index is repeated in the sequence $s_{(1)}, s_{(2)}, \ldots,$ the ordering of corresponding pivot-positions does not affect the final results, such as in Table 6.1. Moreover, if all components in a row are visted the same number of times, then the maximum probability component of the row remains the same.*

## 6.2 Two basic lemmas

Recall from Section 3, that obtaining each term of the sequence $P_{(2)}, P_{(3)}, \ldots$ corresponds to an event of a discrete-event simulation. Our first lemma, Lemma 4 generalizes Table 6.1 above. Our second lemma, Lemma 5, is about the convergence of the probabilities contained in each row of $P_{(t)}$ as $t \to \infty$. We will state Lemma 4 with particular reference to row 1, just to simplify the exposition; one could readily do the same for any other row.

**Lemma** *Assume $P_{(1)}$ has dimension n, and suppose that row 1 of initial matrix $P_{(1)}$ is visited n–1 times inside a subswequence of the sequence of indices $s_{(1)}, s_{(2)}, \ldots$ as follows: position $(1,2)$ is pivot-position for first visit to row 1; position $(1,3)$ is pivot-position for the second visit to row 1; . . . and position $(1,n)$ is pivot-position for the $(n-1)$-th visit to row 1. For j=2,...,n, items (i), (ii) & (iii) below display what row 1 becomes immediately after position $(1,j)$ is processed as a pivot-position.*

*Item (i)–immediately after (1,2) is processed as pivot-position, the components of row 1 become transformed as follows:*

$p_{1,1}$ *is still 0*

$p_{1,2}$ *becomes* $(1+\varepsilon)p_{1,2}/(1+\varepsilon p_{1,2})$

$p_{1,3}$ *becomes* $p_{1,3}/(1+\varepsilon p_{1,2})$

$\vdots$

$p_{1,n}$ *becomes* $p_{1,4}/(1+\varepsilon p_{1,2})$

*Item (ii)–for j=3,...,n-1, immediately after (1,j) is processed as pivot-position, the components of row 1 become transformed as follows:*

$p_{1,1}$ *is still 0*

$p_{1,2}$ *becomes* $(1+\varepsilon)p_{1,2}/(1+\varepsilon p_{1,2}+\ldots+\varepsilon p_{1,j})$

$\vdots$

$p_{1,j}$ *becomes* $(1+\varepsilon)p_{1,j}/(1+\varepsilon p_{1,2}+\ldots+\varepsilon p_{1,j})$

$p_{1,j+1}$ *becomes* $p_{1,j+1}/(1+\varepsilon p_{1,2}+\ldots+\varepsilon p_{1,j})$

$\vdots$

$p_{1,n}$ *becomes* $p_{1,n}/(1+\varepsilon p_{1,2}+\ldots+\varepsilon p_{1,j})$

*Item (iii)–immediately after (1,n) is processed as pivot-position, the components of row 1 become transformed as follows:*

$p_{1,1}$ *is still 0*

$p_{1,2}$ *becomes* $p_{1,2}$ *again*

$\vdots$

$p_{1,n}$ *becomes* $p_{1,n}$ *again*

.

**Proof** This is a straightforward generalization of Table 6.1.

.

**Lemma** *For fixed positive number $\varepsilon \in (0,1)$, define continuous function $f_\varepsilon : [0,1] \to [0,1]$ by $f_\varepsilon(x) = (1+\varepsilon)x/(1+\varepsilon x)$.*

*(a) The numbers 0 & 1 are the only fixed-points of the function $f_\varepsilon$.*

*(b) For any iteration starting number* $x_0 \in (0,1)$, *the fixed-point iteration sequence* $\{x_{k+1} \leftarrow f_\varepsilon(x_k), k = 0,1,\ldots\}$ *converges to 1.*

**Proof** Part (a) is clear from straightforward arithmetic. For Part (b), note that, for every $x \in (0,1)$, we have $x < f_\varepsilon(x) < 1$. Therefore $x_k < f_\varepsilon(x_k) = x_{k+1}$, so that $x_0 < x_1 < \ldots < 1$; and that implies that the sequence $\{x_{k+1} \leftarrow f_\varepsilon(x_k), k = 0,1,\ldots\}$ has 1 as its only limit point.

.

## 6.3 A concluding lemma

In this subsection, in the interest of notation tidiness, we will use a stochastic matrix $Q$ as surrogate for the initial matrix $P_{(1)}$ of Section 3. We will state a concluding lemma that combines conclusions from Lemma 4 and Lemma 5 to demonstrate that: (i) if the simulation procedure (of Section 3) visits a row of matrix $Q$ infinitely many times, then the maximum probability in that row will tend to 1 (as limit) at the end of the simulation run; and (ii) for any starting row index, the computed 1's trace out a single cycle in $Q$'s Markov chain network.

**Lemma** *(a) If pseudo-random variables generated during a simulation run (as indicated in Section 3 above) point-to a component, say component* $j$, *of a row of* $Q$, *say row* $i$, *infinitely many times, then component j of row i has to be the component that contains the maximum probability number in row* $i$, *and that maximum probability converges to 1, at the end of that simulation run.*

*(b) At the end of the simulation run, the computed 1's in matrix Q trace out (that is, describe ) a single cycle in Q's Markov chain network; moreover, that cycle is unique for each simulation-starting row.*

**Proof** For (a): First note that we can choose the simulation parameter ε in Section 3 sufficiently small so as to ensure that the maximum probability in row 1 remains maximal throughout. Therefore, if the simulation pseudo-random variables point-to some components of a row of $Q$ infinitely many times, then the maximal probability component (that is, the component containing the maximum probability number in that row) must be one of such components.

Accordingly, by Lemma 5, if the simulation pseudo-random variables point-to some components of a row of $Q$ infinitely many times, then the maximal probability in that row will converge to 1, at the end of the simulation run, since every row of $Q$ sums up to 1 throughout the simulation run.

For (b), we first note that the pseudo-random numbers (of the simulation run) never point-to a diagonal of $Q$, because the Markov chain represented by $Q$ is a continuous-time Markov chain (so that every diagonal element of $Q$ is 0). Accordingly, $s_{(i+1)} \neq s_{(i)}$ in Section 3, which translates here into saying that the simulation run always changes from each visited row of $Q$ to a *different* row of $Q$. Since the dimension of $Q$ is finite, by following computed 1's in $Q$, one will then go from row to row, and one row will eventually be repeated, and thereby result in a single cycle of states in the network of $Q$ (as each row corresponds to a state).

.

Thus, our algorithm may be described as a calculator that accepts $P_{(1)}$ and $s_{(1)}$ as inputs, and then produces as output a unique cycle among the states of $P_{(1)}$.

# 7. Directions for further work

The main aim of this article is to demonstrate the usefulness of explicitly computing limits of certain practical non-stationary Markov chains. The two broad application instances described in this article are only a tip of the iceberg regarding possible applications. Other application instances include: the maintenance of large energy supply facilities and installations; some stock markets; etc.

In current literature, there are ample predictions that research on various non-stationary Markov chains will become important in the near future, as fast computing will enable practical Markov chain simulation projects to become more prescriptive rather than descriptive.

.

**Acknowledgement** *We hereby express our appreciation and thanks for the indirect contribution of our healthcare project colleagues A. Vanli, C. Blake, and S. Kirtonia when they made constructive comments on aspects of our initial description of the nutritional healhcare center instance.*

# 8. References

.

.

## 9. Appendix - SolveBE functions (3 functions total)

.

```
function S = SolveBE(P) % This is the driver function %
% Pre: a given stochastic matrix P
% Post: a matrix S of solutions of the corresponding BE
[rowP,colP]=size(P);
[Columns,Tr,index]=ForwardPass(P);
S=zeros(index,rowP); % dummy initialization
if index>1 disp('There are several solutions'); end;
%——
for k=1:index % computing solution #k
.   S(k,k)=1; % Recall that S=0, to start
.   for j=(index+1):rowP % Recall that 'Columns' is zero from
. .    u=0; % column #1 through column #index
. .    .for i=1:j-1
```

```
      u=u+(Columns(i,j))*S(k,i); % compute u
    end;
    u=u/(1-Columns(j,j));
  %——
    for i=1:j-1
      S(k,i)=S(k,i)/(1+u); % update old x(j)'s
    end;
  %——
    S(k,j)=u/(1+u); % compute new x(j)
  end;
end;
S=S*Tr;
——
——
End of SolveBE fn, and start of ForwardPass fn
——
——
function [Columns,Tr,index] = ForwardPass(P)
% Pre: P a stochastic matrix
% Post: Columns is matrix of last columns as the forward
% pass is carried out; Tr is matrix indicating re-enumeration
% of states of (markov chain) that is needed to continue
% forward pass; index is size of last identity matrix obtained
% from forward pass (index=1 indicates unique solution of BE)
[rowP,colP]=size(P); Columns=zeros(rowP); Tr=eye(rowP); index=1;
Q=P;
if Q==eye(rowP) % The trivial case
  Columns=0; Tr=0; index=0;
  disp('This is trivially solved');
else
  counter=rowP;
  while index==1 % Repeat until the very end;
    if sum(diag(Q-eye(counter)))==0 break; % Q=1 when counter = 1
    end;
    [Q,T,lastcolumn]=Reduce(Q);
    [rowT,colT]=size(T);
    if rowP>rowT % Other than the start, pad T and lastcolumn
```

```
            T=[T zeros(rowT,rowP-rowT); zeros(rowP-rowT,rowT) eye(rowP-rowT,rowP-rowT)];
            lastcolumn=[lastcolumn;zeros(rowP-rowT,1)];
        end;
        Tr=T*Tr; Columns=T*Columns; % Columns too must be transformed
        Columns(:,counter)=lastcolumn;
        counter=counter-1;
    end;
    index=max(counter,1);
end;
```

——

——

End of ForwardPass fn, and start of Reduce fn

——

——

```
function [Q,T,lastcolumn] = Reduce(P)
% pre: a stochastic matrix P
% post: a matrix Q which is the result of reducing
% P as in Algorithm I, with a matrix T generated
% to reflect re-enumeration as in Algorithm I. The
% inverse of T will be used to transform the resultant
% vector of stationary probabilities as needed. lastcolumn is
% a column vector of last column of P or its re-arrangement
[rowP, colP] = size(P); n=rowP;
T=eye(n); M=T([n 1:(n-1)],[1:n]);
while P(n,n)==1 % change the enumeration of states if necessary
    P=P([n 1:(n-1)],[n 1:(n-1)]); T=M*T;
end;
lastcolumn=P(:,n); % last column of P, with P(rowP,rowP)~=1
Q = zeros(n-1,n-1); % dummy initialization of Q
d=0; % compute 1-p(n,n) indirectly and denote it as d
for i=1:n-1 % so as to aviod doing subtractions
    d=d+P(n,i);
end;
for i = 1:n-1 % update now to obtain Q
    for j = 1:n-1
        Q(i,j) = P(i,j) + P(i,n)*P(n,j)/d;
    end;
end;
```